\documentclass{article}
\usepackage[leqno]{amsmath}
\usepackage{amsfonts}
\begin{document}

\begin{center} {\large \bf Statistics on Multisets} \\ \hfill \\  
Shashikant Mulay and Carl Wagner \\
{\small The University of Tennessee, Knoxville, TN 37996, USA} \end{center} 

\vspace{0.25in}

\noindent {\bf Abstract:} We offer a new proof that a certain $q$-analogue of multinomial 
coefficients furnishes a $q$-counting of the set of permutations of an associated 
multiset of positive integers, according to the number of inversions in such arrangements. 
Our proof uses the fact that such $q$-multinomial coefficients enumerate certain classes of 
chains of subspaces of a finite dimensional vector space over a finite field of cardinality $q$. 
Additionally, we investigate the function that counts the number of permutations of a multiset 
having a fixed number of inversions. \\ 



\noindent {\bf 1. Introduction.} \\

The notational conventions of this paper are as follows: ${\Bbb N}$ and ${\Bbb P}$ denote the set
of non-negative integers and the set of positive integers respectively. If $n \in {\Bbb P}$ and $x$
is an indeterminate, $[n]: = \{1, \dots , n \}$, 
\[n_{x} \, := \, \sum _{i = 0}^{n - 1} x^{i} \;\;\;\;\; \mbox{and} \;\;\;\;\; 
n _{x}^{!} \, := \, \prod _{m = 1}^{n} m_{x}. \]
We follow the convention that  $[0] = \emptyset$, $0_{x} = 0$ and $0_{x}^{!} = 1$. If $n \in {\Bbb N}$,
$r \in {\Bbb P}$ and $(e_{1}, \dots , e_{r}) \in {\Bbb N}^{r}$ is such that $e_{1} + \cdots + e _{r} \leq n$
then define 
\[\tag{1.1} \binom{n}{e_{1}, \dots , e_{r}} _{x} \, := \, 
\frac{n _{x}^{!}}{\prod _{i=1}^{r + 1} (e_{i})_{x}^{!}} \;\;\;\; \mbox{where $e_{r+1} := n - \sum _{i=1}^{r} e_{i}$.} \]
Note that if $r = n - 1$ and $e_{1} = \cdots = e_{r} = 1$, then 
\[\binom{n}{e_{1}, \dots , e_{n - 1}} _{x} \, = \, n_{x}^{!}. \]
When $r =  1$, we write $e_{1} = e$ and so,
\[\tag{1.2}  \binom{n}{e} _{x} \, = \, \frac{n_{x}^{!}}{e_{x}^{!} (n - e)_{x}^{!}} . \]
Since $\binom{n}{0} _{x} = 1 = \binom{n}{n} _{x}$ for all $n \in {\Bbb N}$, and (1.2) implies that
\[\tag{1.3}  \binom{n}{e} _{x} \, = \, \binom{n - 1}{e - 1} _{x} + x^e \binom{n - 1}{e} _{x} 
\;\;\;\;\;\mbox{for $0 < e < n$,} \]
it follows that $\binom{n}{e} _{x}$ is a polynomial in $x$ of degree $e (n - e)$ with coefficients in ${\Bbb N}$. 
Furthermore, letting $e_{0} := 0$, observe that 
\[ \binom{n}{e_{1}, \dots , e_{r}} _{x} \, = \, \prod _{i = 1}^{r} \binom{n - e_{1} - \cdots - e_{i -1}}{e_{i}} _{x}. \]
Hence $\binom{n}{e_{1}, \dots , e_{r}} _{x}$ is also a polynomial in $x$ having coefficients in ${\Bbb N}$ and of degree
\[ \mbox{deg} \binom{n}{e_{1}, \dots , e_{r}} _{x} \, = \, \binom{n}{2} - \sum _{i=1}^{r + 1}  \binom{e_{i}}{2} \, = \,
\sum _{1 \leq i < j \leq r + 1} e_{i} e_{j} . \]
When $x = 1$, the polynomials $n_{x}$, $n_{x}^{!}$, $\binom{n}{e} _{x}$ and $\binom{n}{e_{1}, \dots , e_{r}} _{x}$ evaluate,
respectively, to $n$, $n!$, the binomial coefficient $\binom{n}{e}$, and the $(r + 1)$-nomial coefficient (in abbreviated
notation) $\binom{n}{e_{1}, \dots , e_{r}}$. We will have more to say about this apparently trivial observation in what follows. \\

\noindent {\bf 2. Combinatorial Statistics.} \\

Suppose that $\Delta$ is a finite set of discrete structures and $s: \Delta \rightarrow {\Bbb N}$ is a statistic on $\Delta$
that records some nonnegative integral property of each structure $\delta \in \Delta$. The {\it distribution polynomial} 
$p(\Delta, s, x)$ of $s$ is defined by 
\[\tag{2.1} p(\Delta, s, x) \, := \, \sum _{\delta \in \Delta} x^{s(\delta)} \, = \, 
\sum _{j \in {\Bbb N}} | \{ \delta \in \Delta \, \mid \, s(\delta) = j \} | x^{j} . \]
Of course, $p(\Delta, s, 1) = | \Delta |$ and, if $\Delta$ is equipped with the uniform probability distribution, then
the expected value $\mu_{s}$ of $s$ is equal to $ |\Delta |^{-1} \cdot D_{x} p(\Delta, s, x) | _{x = 1}$. 

As an example, suppose that $n$ and $e$ are positive integers, with $e < n$, and let $M(e, n - e)$ 
denote the family of all multisets containing at most $n - e$ positive integers, each of which is no 
larger than $e$. Note that $|M(e, n - e)| = \binom{n}{e}$. For each $M \in M(e, n - e)$, let $\Sigma (M)$ denote the 
sum of all members of $M$. Then 
\[\tag{2.2}  \binom{n}{e} _{x} \, = \, \sum _{M \in M(e, n - e)} x^{\Sigma (M)} , \]
which may also be expressed in the possibly more familiar form
\[\tag{2.3}  \binom{n}{e} _{x} \, = \, \sum _{m \in {\Bbb N}} p(e, n - e, m) x^{m} , \]
where $ p(e, n - e, m)$ denotes the number of partitions of the integer $m$ with at most $n - e$ parts, each
part being no larger than $e$. A proof of (2.3) appears as early as 1882 in a paper of Sylvester and Franklin
[6, p.269]. More recently Knuth [3] has proved the polynomial identity (2.2) by showing that
\[\tag{2.4}  \binom{n}{e} _{q} \, = \, \sum_{M \in M(e, n - e)} q^{\Sigma (M)} , \]
whenever the integer $q$ is a power of a prime number. Knuth first notes that the $q$-binomial coefficient
appearing on the left in (2.4), enumerates the family ${\cal E}$ of all $e$-dimensional subspaces of an
$n$-dimensional vector space over a finite field of cardinality $q$. He then makes use of the unique $e \times n$
row-echelon matrix associated with each such subspace to define a natural mapping 
$\kappa : {\cal E} \rightarrow M(e, n - e)$ for which each $M \in M(e, n - e)$ has exactly $q^{\Sigma (M)}$ 
preimages with respect to $\kappa$. 

Our aim in this paper is to present Knuth-type proofs for the well-known identity (see Stanley [5, p.21])
\[\tag{2.5} n^{!} _{x} \, = \, \sum_{\theta \in S_{n}} x^{inv(\theta)} , \]
where $S_{n}$ denotes the set of permutations of $[n]$ and for $\theta \in S_{n}$, $inv (\theta)$ denotes the 
number of {\it inversions} of $\theta$, {\it i.e.}, the number of ordered pairs $(i, j)$ with 
$1 \leq i < j \leq n$ such that $\theta (i) > \theta (j)$; and its generalization (see Stanley [5, p.26])
\[\tag{2.6} \binom{n}{e_{1}, \dots , e_{r}} _{x} \, = \, \sum_{\theta \in S(M)} x^{inv(\theta)} , \]
where $S(M)$ is the set of permutations of the multiset 
\[M \, := \, \{1^{e_{1}}, 2^{e_{2}}, \dots , r^{e_{r}}, (r + 1)^{e_{r+1}} \} \]
consisting of $e_{i}$ copies of $i$, with $e_{i} \in {\Bbb P}$, for $1 \leq i \leq r+1$, and 
$e_{1} + e_{2} + \cdots + e_{r} + e_{r+1} = n$. For the above multiset $M$, define $d$ to be the sequence 
$d_{1} < \cdots < d_{r}$, where $d_{i} := e_{1} + \cdots + e_{i}$ for $1 \leq i \leq r$. Our proofs use the fact
that the $q$-factorial appearing on the left in (2.5) and the $q$-multinomial coefficient appearing on the left 
in (2.6) each enumerate chains $V_{1} \subseteq V_{2} \subseteq \cdots \subseteq V_{r}$ of subspaces (so-called 
{\it flags} in the language of algebraic geometry) of an $n$-dimensional vector space over a finite field of 
cardinality $q$ such that $V_{i}$ has dimension $d_{i}$ for $1 \leq i \leq r$. Underlying our proof is a
cell-decomposition of flag spaces. Of course, such cell-decompositions have been known in algebraic geometry 
since the 19'th century, to the best of our knowledge, however, no proof such as the one we provide here has 
appeared in the literature. In the last section we investigate the heretofore neglected function $I_{n}(d; j)$ 
which counts the number of permutations of $M$ having exactly $j$ inversions; this affords a generalization 
of the more familiar function $I_{n}(j)$ which counts the number of permutations of $[n]$ having exactly $j$ 
inversions. The equations, the estimates and the other properties we establish in Theorem 4 merely scratch 
the surface of a potentially deeper combinatorial analysis of $I_{n}(d; j)$. \\

\noindent {\bf 3. Cell decomposition of flag spaces.} \\

Let $n$ be a positive integer. Henceforth we tacitly assume that $n \geq 2$.
Let $V$ be an $n$-dimensional vector space over a (not necessarily finite) field $k$. 
Given a strictly increasing sequence $d : d_1 < d_2 < \dots < d_r$, a sequence 
$V_{1} \subset V_{2} \subset \cdots \subset V_{r}$ of $k$-subspaces of $V$ is called a \emph{$d$-flag} in $V$ 
if $\dim V_i = d_i$ for $1 \leq i \leq r$. The set of all $d$-flags in $V$, called the \emph{space 
of $d$-flags in $V$,} is denoted by $\mbox{FL}(d, V)$.
If $r = 1$, \emph{i.e.}, when the sequence $d$ consists of a single integer (also 
denoted by $d$), the corresponding set $\mbox{FL}(d, V)$ is usually denoted by 
$\mbox{Gr}(d, V)$ and it is called the \emph{Grassmannian} of $d$-dimensional subspaces 
of $V$. If $r = n$, \emph{i.e.}, when $d$ is the sequence $1 < 2 < \dots < n$, 
the corresponding $d$-flag is called a \emph{full flag} in $V$. The space of
full flags in $V$ is usually denoted simply by $\mbox{FL} (V)$. \\

\noindent{\bf Remark:} If $r = 1$ and $d_{1} = n$, then $\mbox{FL}(d, V) = \{V \}$. Note that if $r \geq 2$ 
and $d_{r} = n$, then $\mbox{FL}(d, V)$ is in bijective correspondence with the flag-space 
$\mbox{FL}(d_{1} < \cdots < d_{r - 1}, V)$. So, it suffices to restrict consideration
to the spaces $\mbox{FL}(d, V)$, where $d$ satisfies the additional requirement that $d_{r} < n$. 
By choosing a $k$-isomorphism of vector spaces, $V$ can be identified with $k^{n}$ and then this
induces an identification of $\mbox{FL}(d, V)$ with $\mbox{FL}(d, k^{n})$. In what follows, we 
tacitly assume $V = k^{n}$ and $d_{r} < n$ (whence $r < n)$. \\ 

For positive integers $r, s$ let ${\Bbb M}(r, s, k)$ be the vector space over $k$ of all $r \times s$ 
matrices with entries in $k$. Let ${\Bbb M}(r, k) :=  {\Bbb M}(r, r, k)$, and as usual, let 
$\mbox{GL} (r, k)$ be the multiplicative group of $r \times r$ invertible matrices with entries in $k$.  

Let $d: d_1 < d_2 < \dots < d_r$ be as above. Define $d_{0} := 0$, $d_{r+1} := n$ and 
$e_{i} := d_{i} - d_{i - 1}$ for $1 \leq i \leq r + 1$. Then, $e_{1} + \cdots + e_{r + 1} = n$ and
since $d_{r} < n$, $e_{r + 1} \geq 1$. Define $P(n, d, k)$ to be the set of all $g \in \mbox{GL}(n, k)$ 
such that $g$ is an $(r + 1) \times (r + 1)$ block-matrix $[M_{ij}]$, where $M_{ij} = 0$ for 
$1 \leq j < i \leq r + 1$ and $M_{ii} \in \mbox{GL}(e_{i}, k)$ for $1 \leq i \leq r + 1$. Observe that 
$P(n, d, k)$ is a subgroup of $\mbox{GL}(n, k)$ and $P(n, 1 < 2 < \cdots < n - 1, k)$ is the subgroup 
of upper-triangular matrices. Also, $P(n, 1 < 2 < \cdots < n - 1, k) \, \leq \, P(n, d, k)$ for all 
sequences $d$. 

For a nonnegative integer $e$, let ${\cal R}(n, e, k) \subset {\Bbb M}(n, e, k)$ be the subset of 
matrices of rank $e$. If $A \in {\cal R}(n, e, k),$ then the column-space of $A$, denoted by $C(A)$, is an 
$e$-dimensional $k$-subspace of $k^{n}$. Conversely, any $e$-dimensional $k$-subspace of $k^{n}$ is the 
column-space of some $A \in {\cal R}(n, e, k)$. Furthermore, given $B \in {\cal R}(n, e, k)$, we have
$C(B) = C(A)$ if and only if $A = B g$ for some $g \in \mbox{GL}(e, k)$.  

Fix a sequence $d$ as above and let $e_{1}, \dots , e_{r + 1}$ be the sequence associated to $d$ (as defined above). 
An element of $k^{n}$ is thought of as an $n$-rowed column-matrix. Given $A \in \mbox{GL}(n, k)$, write 
$A := [A_{1}, \dots , A_{r+1}]$ with the understanding that $A_{j}$ is the $n \times e_{j}$ matrix 
made up of columns $d_{j-1} + 1, \dots , d_{j}$ of $A$. By $\Phi (A)$ denote the $d$-flag 
$V_{1} \subset \cdots \subset V_{r}$, where 
\[ V_{m} \, := \, \sum _{i = 1}^{m} C(A_{i}) \;\;\;\;\; \mbox{for $1 \leq i \leq m$.} \]
It is important to note that the sum appearing on the right in the above equation is an internal
direct sum of subspaces of $k^{n}$. \\

\noindent{\bf Theorem 1:} The following holds.
\begin{enumerate}
\item The map $\Phi: \mbox{GL}(n, k) \rightarrow \mbox{FL} (d, V)$ given by $A \rightarrow \Phi (A)$ is surjective.
\item $\Phi (A) = \Phi (B)$ if and only if $A = B g$ for some $g \in P(n, d, k)$. \\
\end{enumerate}

\noindent{\bf Proof:} Fix a $d$-flag ${\cal F}: V_{1} \subset \cdots \subset V_{r}$. Then, find an ordered 
$k$-basis $b(1), \dots , b(n)$ of $k^{n}$ such that $V_{m} = \oplus _{1 \leq i \leq d_{m}} k \cdot b(i)$
for $1 \leq m \leq r$. Define $A_{j + 1}$ to be the $n \times e_{j + 1}$ matrix $[b(d_{j} + 1), \dots , b(d_{j + 1})]$
for $0 \leq j \leq r$ and let $A := [A_{1}, \dots , A_{r+1}]$. Clearly, $A \in  \mbox{GL}(n, k)$ and $\Phi (A) = {\cal F}$. Thus
$\Phi$ is surjective. Suppose $g \in P(n, d, k)$ is the $(r + 1) \times (r + 1)$ block-upper-triangular matrix $[M_{ij}]$
(as in the definition of $P(n, d, k)$ ) and let $[B_{1}, \dots , B_{r+1}] =: B \in \mbox{GL}(n, k)$. Then, 
$B g = [B^{*} _{1}, \dots , B^{*}_{r + 1}]$,
where \[ B^{*} _{j} \, := \, B_{1} M_{1j} + B_{2} M_{2j} + \cdots + B_{j} M_{jj} \;\;\;\; \mbox{for $1 \leq j \leq r + 1$.}\]
Since $M_{jj}$ is invertible, $C(B_{j} M_{jj}) = C(B_{j})$ for $1 \leq j \leq r + 1$. By a straightforward induction,
\[C(B^{*}_{1}) + \cdots + C(B^{*}_{m}) \, = \, C(B_{1}) + \cdots + C(B_{m}) \;\;\;\; \mbox{for $1 \leq m \leq r + 1$.} \]
Hence $\Phi (B) = \Phi (B g)$. Conversely, suppose 
\[A \, := \, [A_{1}, \dots , A_{r+1}] \;\;\;\mbox{and}\;\;\; [B_{1}, \dots , B_{r+1}]\, =: \, B \in \mbox{GL}(n, k) \]
are such that $\Phi (A) = \Phi (B) :=  V_{1} \subset \cdots \subset V_{r}$. Let $V_{r + 1} := k^{n}$. Now since 
$C(A_{1}) = V_{1} = C(B_{1})$, there is a $M_{11} \in \mbox{GL}(e_{1}, k)$ such that $A_{1} = B_{1} M_{11}$. Inductively,
assume that $A_{j} \, = \, B_{1} M_{1j} + \cdots + B_{m} M_{jj}$, where $M_{ij} \in {\Bbb M}(e_{i}, e_{j}, k)$ and 
$M_{jj} \in \mbox{GL}(e_{j}, k)$ for $1 \leq i \leq j \leq m \leq r$. Since $\Phi (A) = \Phi (B)$, we have
\[ V_{m} \oplus C(A_{m + 1}) \; = \; \oplus _{i = 1}^{m + 1} C(A_{i}) \; = \; V_{m + 1} \; = 
\; \oplus _{i = 1}^{m + 1} C(B_{i})\; = \; V_{m} \oplus C(B_{m + 1}) \]
and hence there are natural $k$-linear isomorphisms
\[ C(A_{m + 1})\;  \cong  \; \frac{V_{m + 1}}{V_{m}} \; \cong \; C(B_{m + 1}). \]
Consequently, there exists $M_{(m + 1)(m + 1)} \in \mbox{GL}(e_{m + 1}, k)$ such that 
\[ C(A_{m + 1} - B_{m + 1} M_{(m + 1)(m + 1)}) \, \subseteq \, V_{m}. \] 
In other words, for $1 \leq i \leq m$, there exist matrices $M_{ij} \in {\Bbb M}(e_{i}, e_{j}, k)$ with 
$A_{m + 1} = B_{1} M_{1(m +  1)} + \cdots + B_{m + 1} M_{(m + 1)(m + 1)}$. 
This proves 2. $\Box$ \\  

\noindent{\bf Definitions:} Fix a positive integer $e$ not exceeding $n$. Let $s : s_1 < \cdots < s_e$ be a 
sequence of integers in $[n]$. 
\begin{enumerate}
\item Define $A[s] := [a_{ij}] \in {\cal R}(n, e, k)$ by setting
\[a_{ij}\;=\; \left \{\begin{array}{ll} 1 & \mbox{if $i = s_j,$}
\\ 0 & \mbox{if $i \neq s_j$ } \end{array} \right . \;\;\;\;\; \mbox{for all $(i, j) \in [n] \times [e]$.} \]
\item $M := [m_{ij}] \in {\Bbb M}(n, e, k) $ is said to be in \emph{$s$-reduced form} (resp. {\it anti 
$s$-reduced form}) provided
\[m_{ij}\;=\; \left \{\begin{array}{ll} 1 & \mbox{if $i = s_j$,}
\\ 0 & \mbox{if $i = s_p$ for $p \neq j$ and} \\ 0 & \mbox{if $i < s_j$ (resp. $i > s_{j}$)}
\end{array} \right . \] for all $(i, j) \in [n] \times [e]$.
$M$ is said to be an $s$-reduced form of $N \in {\cal R}(n, e, k)$ if $M$ is in $s$-reduced form
and there exists a matrix $g \in \mbox{GL} (e, k)$ such that $M = N g$. \\ 
\end{enumerate}

\noindent{\bf Lemma 1:} The following holds. 
\begin{description} 
\item[(i)] For $s : s_1 < \cdots < s_e$ in $[n]$, the matrix $A[s]$ is in $s$-reduced form (resp. in 
anti $s$-reduced form). 
\item[(ii)] If $s : s_1 < \cdots < s_e$ and $\sigma : \sigma _{1} < \cdots < \sigma _{e}$ are 
in $[n]$ and $M \in {\cal R}(n, e, k)$, $g \in \mbox{GL} (e, k)$ are such that $M$ is in $s$-reduced 
(resp. anti $s$-reduced) form and $M g$ is in $\sigma$-reduced (resp. anti-reduced) form, then $s = \sigma$ 
and $g = I$ (the identity matrix). 
\item[(iii)] Given $N \in {\cal R}(n, e, k)$, there exists a unique sequence 
$s : s_1 < \cdots < s_e$ in $[n]$ such that $N$ has an $s$-reduced (resp. anti $s$-reduced) form. \\ 
\end{description}

\noindent{\bf Proof:} Assertion (i) clearly holds. Assertions (ii),(iii) are verified by a straightforward 
use of column-reduction to obtain the `reduced column-echelon form'; the `anti' versions of (ii), (iii) can 
be established similarly. Below we present an essential outline of the process (inviting the reader
to formulate its `anti' version).

For $1 \leq j \leq e$, we describe a three-step process. Step 1j: by a suitable permutation of the columns 
$j, \dots , e$, ensure that the $j$-th column has a nonzero entry in some (say) $s_{j}$-th row whereas the rows 
above the $s_{j}$-th row have only $0$ as their entry in columns $j, \dots , e$. Step 2j: multiply the $j$-th column 
by the reciprocal of the entry in the $s_{j}$-th row. Step 3j: subtracting suitable multiples of the $j$-th column from 
each of the remaining columns make sure that the entries appearing in the $s_{j}$-th row and columns other than the 
$j$-th column, are all $0$. To obtain the reduced form of a given $N \in {\cal R}(n, e, k)$, perform the above 
three-step process starting from column $j = 1$ of $N$ and then perform the process for column $j = 2$ on the 
updated matrix, and then perform the process for column $j = 3$ on the updated matrix and so on. $\Box$. \\

\noindent{\bf Definitions:} Let $d : d_1 < d_2 < \dots < d_r$ and $e_{1}, \dots , e_{r + 1}$ be as above.
\begin{enumerate}
\item For a positive integer $e$, let $S[e, n]$ denote the set of all sequences $s : s_1 < \cdots < s_e$ in $[n]$.
An $e$-element subset of $[n]$ is viewed (via the natural ordering of its elements) as a member of $S[e, n]$.
\item Given a subset (possibly empty) $H \subseteq [n]$, a sequence $s \in S[e, n]$ and a matrix $M \in {\cal R}(n, e, k)$,
we say $M$ is in $(s, H)$-reduced (resp. anti $(s, H)$-reduced) form if $M$ is in $s$-reduced (resp. anti $s$-reduced) 
form and for each $i \in H$, the $i$-th row of $M$ is $0$. 
\item Given subsets $\sigma _{1}, \dots , \sigma _{m}$ of $[n]$, we write $\sigma _{1} + \cdots + \sigma _{m} = [n]$ to
mean that $\sigma _{1}, \dots , \sigma _{m}$ form a partition of $[n]$.
\item Let $\pi (d)$ be the subset of $S[e_{1}, n] \times S[e_{2}, n] \times \cdots \times S[e_{r + 1}, n]$ consisting
of $(\sigma _{1}, \dots , \sigma _{r+1})$ such that $\sigma _{1} + \cdots + \sigma _{r + 1} = [n]$.
\item A matrix $A \in G := \mbox{GL} (n, k)$ is said to be in $(\sigma _{1}, \dots , \sigma _{r+1})$-form (resp. anti 
$(\sigma _{1}, \dots , \sigma _{r+1})$-form) provided $(\sigma _{1}, \dots , \sigma _{r+1}) \in \pi (d)$ and 
$A = [A_{1}, \dots , A_{r+1}]$, where $A_{j}$ is in $(\sigma _{j}, \; \sigma _{1} \cup \cdots \cup \sigma _{j - 1})$-reduced 
form (resp. anti $(\sigma _{1}, \dots , \sigma _{r+1})$-form) for $1 \leq j \leq r + 1$ (convention: if
$j =1$, then $\sigma _{1} \cup \cdots \cup \sigma _{j - 1} := \emptyset$). \\
\end{enumerate}

\noindent{\bf Remark:} Note that for there to be an $M$ in $(s, H)$-reduced (resp. anti $(s, H)$-reduced) form, it 
is necessary that $s \cap H = \emptyset$. \\

\noindent{\bf Lemma 2:} Let $A \in \mbox{GL} (n, k)$. Then, there exists a unique 
$(\sigma _{1}, \dots , \sigma _{r+1}) \in \pi (d)$ and a unique $g \in P(n, d, k)$ such that
$A g$ is in $(\sigma _{1}, \dots , \sigma _{r+1})$-form (resp. anti $(\sigma _{1}, \dots , \sigma _{r+1})$-form). \\

\noindent{\bf Proof:} Again, the `anti' version of the proof is left to the reader; it is easily obtained by suitable
modification of the following arguments. By induction on $m \leq r + 1$ we find matrices $M_{ij} \in {\Bbb M}(e_{i}, e_{j}, k)$, 
with $M_{ii} \in \mbox{GL} (e_{i}, k)$ such that $B _{j} := A_{1} M_{1j} + \cdots + A_{j} M_{jj}$
is in $(\sigma _{j}, \; \sigma _{1} \cup \cdots \cup \sigma _{j - 1})$-reduced form for 
$1 \leq i \leq j \leq m$. Case $m = 1$: thanks to Lemma 1, there is a unique $M_{11} \in \mbox{GL} (e_{1}, k)$ 
and a unique $\sigma _{1} \in S[e_{1}, n]$ such that $B_{1}:= A_{1} M_{11}$ is in $(\sigma _{1}, \; \emptyset)$-reduced 
form. Case $m \geq 2$: By induction we assume that we have found the desired $M_{ij}$ for $1 \leq i \leq j \leq m - 1$.
For $1 \leq t \leq m - 1$, let $N_{t} \in {\Bbb M}(e_{t}, e_{m}, k)$ be such that for each 
$i \in \sigma _{1} \cup \cdots \cup \sigma _{m - 1} $, the $i$-th row of $C_{m} := A_{m} + (B_{1} N_{1} + \cdots + B_{m - 1} N_{m - 1})$ 
is a zero row. Let $M_{mm} \in \mbox{GL} (e_{m}, k)$ be the unique matrix and let $\sigma _{m} \in S[e_{m}, n]$ be the unique
sequence such that $C_{m} M_{mm}$ is in $\sigma _{m}$-reduced form. Then, $C_{m} M_{mm}$ is automatically
in $(\sigma _{m}, \; \sigma _{1} \cup \cdots \cup \sigma _{m - 1})$-reduced form. Letting $B_{m} := C_{m} M_{mm}$
\[ M_{im} \; :=  \left ( \sum _{q = i} ^{m - 1} M_{iq} N_{q} \right ) M_{mm} \;\;\;\;\mbox{for $1 \leq i \leq m - 1$} \]
we infer that $B_{m} =  A_{1} M_{1m} + \cdots + A_{m} M_{mm}$. Also, letting $g$ be the $(r + 1) \times (r + 1)$ 
block-upper triangular matrix $[M_{ij}]$, we have  $g \in P(n, d, k)$ and 
$[B_{1}, \dots , B_{r + 1}] = [A_{1}, \dots , A_{r + 1}] g$. This establishes the existence part of our assertion. 

We proceed to prove the asserted uniqueness. Suppose $A := [A_{1}, \dots , A_{r + 1}]$ is in 
$(\sigma _{1}, \dots , \sigma _{r+1})$-form and there is $g \in P(n, d, k)$ such that $[B_{1}, \dots , B_{r + 1}] := A g$ 
is in $(\tau _{1}, \dots , \tau _{r+1})$-form. Say $g$ is the $(r + 1) \times (r + 1)$ block-upper triangular matrix $[M_{ij}]$.
By (ii) of Lemma 1, we at once infer that $\sigma _{1} = \tau _{1}$ and $M_{11} = I$. Hence $B_{1} = A_{1}$. By induction,
assume that $\sigma _{j} = \tau _{j}$, $M_{jj} = I$ and $M_{ij} = 0$ for $1 \leq i < j \leq m - 1$. Then, we must also have
$B_{j} = A_{j}$ for $1 \leq j \leq m - 1$ and hence 
\[ B_{m} - A_{m} M_{mm} \; =  \; B_{m - 1} M_{(m - 1)m} + \cdots + B_{1} M_{1m} . \]
Consider a $1 \leq j \leq m - 1$. For each $i \in \sigma _{1} \cup \cdots \cup \sigma _{j - 1} = \tau _{1} \cup \cdots \cup \tau _{j - 1}$, 
the $i$-th row of the matrix on the left (in the above equation) as well as the matrix $B_{m - 1} M_{(m - 1)m} + \cdots + B_{j} M_{jm}$,
is $0$. When $j = 2$, using the fact that $B_{1}$ is $\tau_{1} =\sigma _{1}$-reduced, we get $M_{1m} = 0$. Repeating this argument
for each of $j = 3, \dots, m - 1$, we infer that  $M_{im} = 0$ for $1 \leq i \leq m - 1$. Consequently, $B_{m} = A_{m} M_{mm}$.
Now by (ii) of Lemma 1, $\sigma _{m} = \tau _{m}$ and $M_{mm} = I$. It follows that $g = I$. $\Box$ \\

\noindent{\bf Definition:} Let $d: d_1 < d_2 < \dots < d_r$ and $\pi (d)$ be as above. To  each partition
\[ \sigma  \, := \,  (\sigma _{1}, \dots , \sigma _{r+1}) \,  \in  \, \pi (d) \] 
we associate subsets, or {\it cells},  $W_{\sigma}$ and $ \widehat{W}_{\sigma}$ defined by
\[ \begin{array}{lll} W_{\sigma} & := & \{A \in \mbox{GL} (n, k) \, \mid \, \mbox{$A$ is in $(\sigma _{1}, \dots , \sigma _{r+1})$-form } \} ,
\\  \widehat{W}_{\sigma} & := & \{A \in \mbox{GL} (n, k) \, \mid \, \mbox{$A$ is in anti $(\sigma _{1}, \dots , \sigma _{r+1})$-form } \}.
\end{array} \]

\noindent{\bf Theorem 2:} Given $d: d_1 < d_2 < \dots < d_r$ in $S[r, n]$, the following holds.
\begin{description}
\item[(i)] We have the cell-decompositions
\[\frac{\mbox{GL} (n, k)}{ P(n, d, k)} \;\; \cong  \; \bigsqcup _{\sigma \in \pi (d)} W_{\sigma} \; = \;
\bigsqcup _{\sigma \in \pi (d)} \widehat{W}_{\sigma} . \]
\item[(ii)] Let $\Phi : \mbox{GL} (n, k)  \rightarrow \mbox{FL}(d, V)$ be the map as in Theorem 1. Then,
\[ \mbox{FL}(d, V) \;\; = \; \bigsqcup _{\sigma \in \pi (d)} \left \{ \Phi (A) \, \mid \, A \in W_{\sigma} \right \} 
\; = \; \bigsqcup _{\sigma \in \pi (d)} \left \{ \Phi (A) \, \mid \, A \in \widehat{W}_{\sigma} 
 \right \}. \]
\end{description}

\noindent{\bf Proof:} By Lemma 2, given a left-coset $L$ of $P(n, d, k)$ in $\mbox{GL} (n, k)$, there is a unique 
$\sigma \in \pi (d)$ and a unique $A \in W_{\sigma}$ (resp. $A \in \widehat{W}_{\sigma}$)
such that $L = A P(n, d, k)$. Mapping a left-coset $L$ to its representative yields the bijective correspondence
asserted in (i). In view of Theorem 1, (ii) follows from (i). $\Box$ \\ 

\noindent {\bf 4. Dimension Counting.} \\

Fix $d: d_1 < d_2 < \dots < d_r$ in $S[r, n]$ with $d_{r} < n$. As before, let $d_{0} := 0$ and $d_{r + 1} := n$. 
For notational simplicity, a sequence $s \in S[e, n]$ is henceforth written as $s(1) < s(2) < \cdots < s(e)$.  \\

\noindent{\bf Definitions:} Let $\sigma := (\sigma _{1}, \dots , \sigma _{r+1}) \in \pi (d)$. As before, let $S(M)$ 
be the set of permutations of the multiset $M := \{1^{e_{1}}, 2^{e_{2}}, \dots , r^{e_{r}}, (r + 1)^{e_{r+1}} \}$. 
\begin{enumerate}
\item For $j \in [n]$, define 
\[ \sigma (j) \; := \; \sigma _{m} (j - d_{m}) \;\;\;\;\mbox{if $d_{m} < j \leq d_{m + 1}$ for some $0 \leq m \leq r$.} \]
\item For $j \in [n]$, let $\mu (j) := m + 1$ provided $d_{m} < j \leq d_{m + 1}$ with $0 \leq m \leq r$.
\item For $0 \leq m \leq r + 1$, let $T (m) := [n] \setminus \{\sigma (i) \, \mid \, 1 \leq i \leq d_{m} \}$.
\item For $j \in [n]$, let $\delta (\sigma, j) := |\Delta (\sigma, j)|$, where
 \[\Delta (\sigma, j) \, := \, \{ t \in T(\mu (j)) \, \mid \, t > \sigma (j) \}. \] 
\item Let $\lambda (\sigma) := \sum _{j = 1}^{j = n} \delta (\sigma, j)$. 
\item Define $\Theta _{\sigma} \in S(M)$ by $\Theta _{\sigma} (i) := n - \sigma (i) + 1$ for all $i \in [n]$. \\
\end{enumerate} 

\noindent{\bf Remarks:} Let $\nu (d) := \sum _{1 \leq i < j \leq r + 1} e_{i} e_{j}$.
\begin{enumerate} 
\item Observe that $\sigma$, as defined in the first of the above definitions, is a permutation
of $[n]$. The corresponding permutation $\Theta _{\sigma}$ is also called the `opposite' or the `dual' of $\sigma$. 
\item Clearly, if $m := \mu (j)$, then $\delta (\sigma , j) \leq |T(m)| = n - d_{m} = e_{m+1} + \cdots + e_{r+1}$. 
In particular, $\delta (\sigma , j) = 0$ for $j > d_{r}$ and $\lambda (\sigma) \leq \nu (d)$. Moreover,
$\delta (\sigma , j) = n - d_{\mu (j)}$ for $1 \leq j \leq d_{r}$ if and only if $\lambda (\sigma) = \nu (d)$ if and 
only if $\sigma$ is the identity permutation. At the opposite extreme, $\delta (\sigma , i) = 0$ for $1 \leq i \leq j$ 
if and only if $\lambda (\sigma) = 0$ if and only if $\sigma (i): n - d_{i} + 1 < \cdots < n - d_{i-1}$ 
for $1 \leq i \leq r + 1$.
\item Consider the case where $r = 1$, {\it i.e.}, $e_{1} = d_{1}$ and $e_{2} = n - d_{1}$. Suppose $k$ is an integer 
with $0 \leq k \leq \nu (d) = e_{1}e_{2}$. Let $a, b \in {\Bbb N}$ be such that $k = a e_{2} + b$ and $b < e_{2}$. 
Let $\tau:= (\tau _{1}, \tau _{2}) \in \pi (d)$ be such that  
\[ \tau _{1} \, = \, \{j \, \mid \, 1 \leq j \leq a \} \cup \{n - j \, \mid \, 0 \leq j \leq e_{1} - a - 1 \} \cup 
\{ n - e_{1} + a - b \} . \]
Then, it is easy to verify that $\lambda (\tau) = k$. \\
\end{enumerate}

\noindent{\bf Theorem 3:} The following holds.
\begin{description} 
\item[(i)] Let $A := [a_{ij}] \in \mbox{GL} (n, k)$ and $\sigma \in \pi (d)$. Then, $A \in W_{\sigma}$ if and only if 
\[a _{ij} \, = \, \left \{ \begin{array}{ll} 1 & \mbox{if $i = \sigma (j)$, } \\ 
0 & \mbox{if $i \neq \sigma (j)$ and $i \not \in \Delta (\sigma, j)$.} \end{array} \right . \]
\item[(ii)] $\Theta : \pi(d) \rightarrow S(M)$ defined by $\sigma \rightarrow \Theta _{\sigma}$ is a bijective map. Moreover,
letting $\theta := \Theta _{\sigma}$, we have $\lambda (\sigma) = inv (\theta)$.
\item[(iii)] Letting $\theta := \Theta _{\sigma}$, we have $\lambda (\sigma) = inv (\theta)$. So,
\[ W_{\sigma} \, \cong \, k^{\lambda (\sigma)} \, = \,  k^{inv (\theta)} . \]
\end{description}

\noindent{\bf Proof:} If $A \in W_{\sigma}$, then it is straightforward to verify that the entries of $A$ 
satisfy (i). Conversely, suppose entries of $A$ satisfy (i). As before,
let $A := [A_{1}, \dots , A_{r + 1}]$ and let $\sigma := (\sigma _{1}, \dots , \sigma _{r+1})$. Let $m$ be
an integer such that $0 \leq m \leq r$. If $d _{m} < j \leq d_{m + 1}$, then the $j$-th column of $A$ is
the $(j - d_{m})$-th column of $A_{m + 1}$. Consider the $q$-th column of $A_{m + 1}$. Letting $j := q + d_{m}$, 
we have $a_{ij} = 1$ if $i = \sigma (j) = \sigma _{m + 1} (q)$, $a_{ij} = 0$ if $i < \sigma (j)= \sigma _{m + 1} (q)$
and $a_{ij} = 0$ if $\sigma (j) = \sigma _{m + 1} (q) < i \in \sigma _{1} \cup \cdots \cup \sigma _{m + 1}$ .
Thus $A_{m + 1}$ is in $(\sigma _{m+1}, \; \sigma _{1} \cup \cdots \cup \sigma _{m})$-reduced form. This
proves (i). 

It is straightforward to verify that $\Theta$ is a bijective map. Fix a $\sigma \in \pi(d)$ and let 
$\theta := \Theta _{\sigma}$. For $1 \leq i < j \leq n$, we have $\mu (i) \leq \mu (j)$ and 
$\theta (i) > \theta (j)$ if and only if $\sigma (i) < \sigma (j)$ if and only if $\sigma (j) \in \Delta(\sigma, i)$.
In other words, for $i \in [n]$, the set $\{j \in [n] \, \mid \, i < j, \;\;\; \theta (i) > \theta (j)\}$ is in 
one-to-one correspondence with the set $\Delta (\sigma, i)$. Now it readily follows that $\lambda (\sigma) = inv (\theta)$.
In view of Theorem 2, (iii) follows from (i) and (ii). $\Box$ \\

\noindent{\bf Corollary:} Let $d: d_{1} < \cdots < d_{r}$ be in $S[r, n]$ with $d_{r} < n$, $d_{0} := 0$, 
$d_{r + 1} := n$ and let $e_{i} := d_{i} - d_{i - 1}$ for $1 \leq i \leq r + 1$. Assume $k$ is a finite field 
with $|k| = q$. Then, the following holds.
\begin{description}
\item[(i)]
\[ |\mbox{FL}(d, V)| \,= \,
\frac{ \prod _{i = 0}^{n - 1} (q^{n} - q^{i})}{\prod _{i = 1}^{r + 1} \prod _{j = 0}^{e_{i} - 1}(q^{e_{i}} - q^{j})
\prod _{1 \leq i < j \leq r + 1 } q^{e_{i} e_{j}} } \, = \, \binom{n}{e_{1}, \dots , e_{r}} _{q} . \] 
\item[(ii)] \[ |\mbox{FL}(d, V)| \; = \; \sum _{\theta \in S(M)} q^{inv(\theta)} .\] 
\item[(iii)] \[ \binom{n}{e_{1}, \dots , e_{r}} _{x} \, = \, \sum _{\theta \in S(M)} x^{inv(\theta)} . \] \\
\end{description}

\noindent{\bf Proof:} By Theorem 2, $|\mbox{FL}(d, V)| = |\mbox{GL} (n, k) / P(n, d, k)|$. Since we have
$|\mbox{GL} (n, k)| = \prod _{i = 0}^{n - 1} (q^{n} - q^{i})$ and
\[ |P(n, d, k)| \, = \,
\left ( \prod _{i = 1}^{r + 1} \prod _{j = 0}^{e_{i} - 1}(q^{e_{i}} - q^{j}) \right )
\left ( \prod _{1 \leq i < j \leq r + 1 } q^{e_{i} e_{j}} \right ) , \]
the first equality in (i) follows. The second equality asserted in (i) is essentially the equality (1.1). 
In view of Theorem 2 and (iii) of Theorem 3, (ii) holds. Since (i) and (ii) hold for infinitely many $q$,
assertion (iii) must hold. $\Box$ \\

\noindent{\bf Remarks:}
\begin{enumerate}
\item It can be easily verified that $\widehat{W}_{\sigma}\, \cong \, k^{inv (\sigma)}$. So, 
decomposing $\mbox{FL}(d, V)$ into cells $\widehat{W}_{\sigma}$, we can identify $\pi(d)$ with 
$S(M)$ in a straightforward manner ({\it i.e.}, without $\Theta$). Our preference for the cells $\widehat{W}_{\sigma}$
is rooted in the belief that $\sigma$-forms (of matrices) are more familiar than anti $\sigma$-forms. 
\item Consider indeterminates $X$ and $z_{ij}$ for $1 \leq j \leq e_{i}$, $1 \leq i \leq r + 1$. 
Let $g_{i} := X^{e_{i}} + \sum _{1 \leq j \leq e_{i}} z_{ij} X^{e_{i} - j}$ for $1 \leq i \leq r + 1$ and let
$g$ be the product of $g_{1}, g_{2}, \dots , g_{r+1}$. For $m \in [n]$, let $a_{m}$ denote the coefficient of 
$X^{n-m}$ in $g$. Let $R$ be the polynomial ring over ${\Bbb Z}$ in the $n$ indeterminates $z_{ij}$ and let 
$J$ denote the ideal of $R$ generated by $a_{1}, \dots , a_{n}$. We let $z_{ij}$ have weight $j$ for all $i, j$.
Then, $J$ is a weighted homogeneous ideal of the weighted homogeneous ring $R$. Furthermore, $J$ is an ideal 
theoretic complete intersection. Then, it is known that the Hilbert series of the weighted (or {\it graded}) ring 
$R / J$ is the polynomial appearing in (iii) of the above corollary. \\
\end{enumerate}

\noindent {\bf 5. Multiset-permutations with fixed number of inversions.} \\

Let $t$ be an indeterminate and let $w:= (w_{1}, \dots , w_{n})$ be an $n$-tuple of positive integers (where $n$ is also 
a positive integer). For $m \in {\Bbb Z}$, let $D_{w}(m)$ be defined by the equation 
\[ \prod _{i = 1}^{n} \frac{1}{(1 - t^{w_{i}})} \, = \, \sum _{m \in {\Bbb Z}} D_{w}(m) t^{m} . \]
Since the rational function on the left is a power-series in $t$ with coefficients in ${\Bbb N}$, we have $D_{w} (m) \in {\Bbb N}$
for all $m \in {\Bbb Z}$ and $D_{w}(m) = 0$ for $m < 0$. Also, observe that $D_{w} (0) = 1$. For $m \in {\Bbb N}$, the integer
$D_{w}(m)$ is known as the {\it Sylvester's denumerant}; clearly,
\[ D_{w} (m) \, = \, |\{(i_{1} , \dots , i_{n}) \in {\Bbb N}^{n} \, \mid \, i_{1} w_{1} + \cdots  + i_{n} w_{n} = m \}| .\]
It is well known that if $\lambda := lcm (w_{1}, \dots , w_{n})$, then for each $j$ with $0 \leq j \leq \lambda - 1$, there is a 
polynomial $Q_{j}(t) \in {\Bbb Q}[t]$ of degree $n - 1$ such that $D_{w}(m) = Q_{j}(m)$ provided $m \equiv j \; mod \; \lambda$. 
So, \[ D_{w}(m) = P_{0}(w; m) + P_{1}(w; m) m + \cdots + P_{n-1}(w; m) m^{n-1}, \] where each $P_{i}(w; m)$ is a ${\Bbb Q}$-valued 
periodic function of $m$. For more on this topic the reader is referred to [1], [2] and their list of references. \\

\noindent{\bf Definitions:} Let $n$ be a positive integer and let $d: d_{1} < \cdots < d_{r}$ be a sequence of positive 
integers with $d_{r} < n$; as before, $d_{0} := 0$, $d_{r + 1} := n$
\begin{enumerate}
\item Given $a, b \in {\Bbb Z}$, let
\[ \binom{a}{b} \, := \, \left \{\begin{array}{ll} 0 & \mbox{if min$\{a, b\} < 0$,} \\ 0 & \mbox{if $a < b$ and} \\
\prod_{0 \leq i \leq b - 1} \frac{(a - i)}{i + 1 } & \mbox{if $0 \leq b \leq a$. } \end{array} \right . \]
\item Let $1_{n}$ denote the $n$-tuple $(w_{1}, \dots , w_{n})$, where $w_{i} = 1$ for $1 \leq i \leq n$. 
\item Let $([n])$ denote the $n$-tuple $(w_{1}, \dots , w_{n})$, where $w_{j} = j$ for $1 \leq j \leq n$. 
\item For a subset $T \subseteq [n]$, let $\omega(T) \, := \, \sum _{i \in T} i$. For $r \in {\Bbb N}$, let
\[\psi_{n}(r) \, := \, \sum _{T \in \Omega (n, r)} (-1)^{|T|} , \;\;\;\mbox{where}\;\;\;
\Omega (n, r) \, := \, \{ T \subseteq [n] \, \mid \, \omega (T) = r \}. \] 
\item Define 
\[ \varepsilon (d, n) \, := \, (w_{1}, \dots , w_{n}),  \;\;\;\; \mbox{where $w_{i} = i - d_{j}$ provided 
$d_{j} < i \leq d_{j+1}$.}  \] 
\item For $k \in {\Bbb N}$, let $I_{n}(d; k)$ denote the number of permutations $\theta \in S(M)$ with 
$inv (\theta) = k$. In the special case where $r = n -1$ (and hence $d: 1 < 2 < \cdots < n - 1$), we write 
$I_{n}(k)$ instead of $I_{k}(d; n)$. It is convenient to allow $d$ to be the empty sequence (in which case,
$\pi (d)$ is (by convention) the trivial subgroup of $S_{n}$). 
\item A sequence $d^{*}: d_{1}^{*} < \cdots < d_{s}^{*} < n$ is called a {\it refinement of $d$} if each $d_{i}$
is a member of $d^{*}$. 
\item For $n \in {\Bbb N}$, let $f_{n}(t):= (1 - t) (1- t^2) \cdots (1 - t^n) \in {\Bbb Z}[t]$ and let
$G_{n}(t) \in {\Bbb Z}[t]$ be defined by $G_{n} (t) := f_{n}(t) / (1 - t)^{n}$. \\
\end{enumerate}

\noindent{\bf Remarks:} 
\begin{enumerate}
\item It is easily seen that $\psi _{n} (r)$ is the coefficient of $t^{r}$ in $f_{n}(t)$. From the 
identity $f_{n}(t) = (-1)^{n} \cdot t^{n (n + 1) /2} \cdot f_{n} (1 / t)$, we deduce that 
\[ \psi _{n}(r) \, = \, (- 1)^{n} \cdot \psi _{n} \left (\frac{n (n + 1)}{2} - r \right ). \]
Applying Euler's pentagonal number theorem, it can be easily verified that
\[\mbox{for $1 \leq r \leq n$,} \;\;\;\; \psi _{n}(r) \, = \, \left \{ 
\begin{array}{ll} (-1)^{s} & \mbox{if $ 2 r = s (3 s \pm 1)$ with $s \in {\Bbb N}$,} \\ \\
0 & \mbox{otherwise.} \end{array} \right . \]
For a positive integer $k$, define the {\it restricted divisor-sum}
\[\sigma _{n} (k) \, := \, \sum _{1 \leq d \leq n,\; d | k} d \, = \, 
\sum _{d = 1}^{ {\rm min} \{n, k \}} \left \lfloor 1 + \left \lfloor \frac{k}{d} \right  \rfloor - \frac{k}{d} \right \rfloor d \]
and let $\alpha_{n} (k) := \sigma _{n} (k) / k$. Then, the (formal) identity $f_{n}(t) = exp(log(f_{n} (t))$ 
provides the formula
\[\psi _{n} (r) \, = \, \sum _{i_{1} + 2 i_{2} + \cdots + m i_{m} = r } \frac{(- 1)^{i_{1} + i_{2} + \cdots + i_{m}}}
{i_{1}! i_{2}! \cdots i_{m}!}\, \alpha_{n}(1)^{i_{1}} \alpha_{n}(2)^{i_{2}} \cdots \alpha_{n}(m)^{i_{m}}. \]
There is yet another such formula that can be derived from the pentagonal number expansion; but it is 
equally complicated and perhaps of little use. We also have the obvious inequality 
\[ | \psi _{n} (r) | \, \leq \, \binom{n - 1 + r}{n - 1} \;\;\;\;\mbox{for all $r \in {\Bbb N}$.} \] 
For $n < r < n (n -1) / 2$, very little seems to be known regarding the size or sign of $\psi _{n} (r)$.
\item Note that $I_{n}(k)$ is the number of permutations of $\{1, \dots , n \}$ having exactly $k$ inversions and so, 
\[\sum _{k \in {\Bbb N}} I_{n}(k) t^{k} \, = \, \binom{n}{1, \dots , 1} _{t} \, = \, G_{n}(t). \]
Since $G_{n}(t) = t^{n (n - 1) /2} \cdot G _{n} (1 / t) = G_{n-1}(t) (1 + t + \cdots + t^{n-1})$, 
\[I_{n}(k) \, = \,I_{n} \left (\frac{n (n - 1)}{2} - k \right ) \;\;\;\; \mbox{and} \;\;\;\; 
 I_{n}(k) \, = \, \sum_{j = {\rm max} \{0, k - n + 1 \}}^{k} I_{n-1}(j). \]
In particular, $I_{n} (k) \geq 1$ for $0 \leq k \leq n (n - 1) /2$. \\ 
\end{enumerate}

\noindent{\bf Theorem 4:} Let the notation be as above; in particular, let $w:= (w_{1}, \dots , w_{n})$ be an $n$-tuple of 
positive integers and let $d: d_{1} < \cdots < d_{r}$ be a sequence of positive integers with $d_{r} < n$. As before,
let $e_{i} := d_{i} - d_{i-1}$ for $1 \leq i \leq r +1$, where $d_{0} = 0$ and $d_{r+1} = n$.
\begin{description}
\item[(i)] For $r \in {\Bbb N}$, we have
\[ \frac{ \prod _{i = 1}^{r} (1 - t^{i}) }{\prod _{i = 1 }^{n} (1 - t^{w_{i}}) }\, = \, 
\sum _{m \in {\Bbb N}} \left ( \sum _{T \subseteq [r]} (-1)^{|T|} D_{w}(m - \omega(T)) \right ) t^{m}  . \]
\item[(ii)] For $k \in {\Bbb N}$, we have
\[I_{n}(d; k) \, = \, \sum _{T \subseteq [n]} (-1)^{|T|} D_{\varepsilon(d, n)}(k - \omega(T)) \, = \,
\sum _{i = 0}^{k} \psi_{n}(i) D_{\varepsilon(d, n)}(k - i) .\]
\item[(iii)] For $k \in {\Bbb N}$, we have
\[I_{n}(k) \, = \, \sum _{T \subseteq [n]} (-1)^{|T|} \binom{n - 1 + k - \omega(T)}{n - 1} \, = \, 
\sum _{i = 0}^{k} \psi_{n}(i) \binom{n - 1 + k - i}{n - 1} . \]  
\item[(iv)] Let $d^{*}$ be a refinement of $d$ and let
\[ {\cal A}_{m} := \{j:=(j_{1}, \dots , j_{r+1}) \in {\Bbb N}^{r+1} \, \mid  \, j_{1} + \cdots + j_{r+1} = m \}. \]
Then, for $1 \leq i \leq r + 1$, there is a (possibly empty) sequence $d[i]^{*}$ of positive integers 
$ < e_{i}$ determined by $d^{*}$ (in a canonical manner). Furthermore,
\[I_{n}(d; k) \, = \, I_{n}(d^{*}; k) - \sum _{m = 1}^{k} \left ( \sum _{j \in {\cal A}_{m}}
\prod _{m = 1}^{r + 1} I_{e_{m}}(d[m]^{*}; j_{m}) \right ) I_{n}(d; k - m)  \] 
for all $k \in {\Bbb N}$. In particular, $I_{n}(d; k) \leq I_{n}(d^{*}; k)$.
\item[(v)] Let $\nu (d) := \sum _{1 \leq i < j \leq r + 1} e_{i} e_{j}$. Then, $\nu (d) \leq n(n-1)/2$,
\[I_{n}(d; k) \, = \, I_{n}(d; \nu (d) - k)\;\;\;\;\; \mbox{for all $k \in {\Bbb N}$} \] 
and $I_{n}(d; k) \geq 1$ for $0 \leq k \leq \nu(d)$. \\
\end{description}

\noindent{\bf Proof:} We prove (i) by induction on $r$. If $r = 0$, then since $[r] = \emptyset$ and 
$\prod _{i = 1}^{r} (1 - t^{i}) = 1$, assertion (i) trivially holds. Fix a positive integer $r$ such that (i)
holds for $r - 1$. If $f := \sum c(m) t^{m} \in {\Bbb Q}[[t]]$, then observe that the coefficient of $t^{m}$ 
in the product $(1 - t^{r}) f$ is $c(m) - c(m - r)$. In particular, the coefficient of $t^{m}$ in the product
\[(1 - t^{k}) \cdot \sum _{m = 0} \left ( \sum _{T \subseteq [r - 1]} (-1)^{|T|} D_{w}(m - \omega(T)) \right ) t^{m} \]
is the difference
\[\sum _{S \subseteq [r - 1]} (-1)^{|S|} D_{w}(m - \omega(S)) - 
\sum _{S \subseteq [r - 1]} (-1)^{|S|} D_{w}(m - \omega(S) - r) .\]
Given a subset $T \subseteq [r]$, letting $S := T \cap [r-1]$, we have either $T = S$ or $T = S \cup \{r \}$.
In the first case $\omega(T) = \omega (S)$ and in the second case, $\omega(T) = \omega(S) + r$. So, in view of our
induction hypothesis, (i) holds for $r$. 

Let $d$ be as in (ii) and let $d_{0} := 0$, $d_{r + 1} := n$ and let $e_{i} := d_{i} - d_{i - 1}$ for 
$1 \leq i \leq r + 1$. Combining our earlier observations, we get
\[ \sum _{k \geq 0} I_{n}(d; k) t^{k} \, = \, \binom{n}{e_{1}, \dots , e_{r}} _{t} \, = \, 
\frac{ \prod _{i = 1}^{n} (1 -  t^{i})}{\prod _{i = 1}^{r + 1} \prod _{j = 1}^{e_{i}}(1 - t^{j})}. \]
Consequently, (ii) is a special case of (i) in which $w = \varepsilon(d, n)$. As is well known, 
\[ D_{1_{n}}(m) \, = \, \binom{n - 1 + m }{n - 1} \;\;\;\;\;\; \mbox{for all $m \in {\Bbb N}$.} \]
Hence (iii) is a special case of (i) in which $w = 1_{n}$.  

As in (iv), let $d^{*}$ be a refinement of $d$. Given $1 \leq i \leq r + 1$, suppose 
$d_{i-1} < d_{j}^{*} < \cdots < d_{j+s}^{*} < d_{i}$ and then, let
\[ d[i]^{*}:\; (d_{j}^{*} - d_{i-1}) < (d_{j+1}^{*} - d_{i-1}) < \cdots < (d_{j+s}^{*} - d_{i-1}) .\]
Note that $d[i]^{*}$ may be empty. Let $e_{1}^{*}, \dots , e_{p+1}^{*}$ be the sequence of
the consecutive differences of the members of $d^{*}$ (with $d_{0}^{*} = 0$ and $d_{p+1}^{*} = n$).
Likewise, for $1 \leq i \leq r + 1$, let $e[i]_{1}^{*}, \dots , e[i]_{p[i] + 1}^{*}$ be the sequence of
the consecutive differences of the members of $d[i]^{*}$ (with $d[i]_{0}^{*} = 0$ and $d[i]_{p[i]+1}^{*} = e_{i}$).
Observe that 
\[\binom{n}{e_{1}, \dots , e_{r+1}} _{t}  \, = \,\binom{n}{e_{1}, \dots , e_{p+1}} _{t} 
\prod _{i = 1}^{r+1} \binom{n}{e[i]_{1}^{*}, \dots , e[i]_{p[i]+1}^{*}} _{t} . \] 
Now (iv) follows by equating the coefficients of like powers on the both sides of this equation. 

Clearly, we have 
\[ \nu (d) \, = \, \frac{n (n + 1)}{2} - \sum _{i = 1}^{r + 1} \frac{e_{i} (e_{i} + 1)}{2} \, = \,
\frac{n (n + 1)}{2} - n + \sum _{i = 1}^{r + 1} \frac{e_{i} (e_{i} - 1)}{2} \, \leq \, \frac{n (n - 1)}{2}. \]
Also, it is straightforward to verify that
\[\binom{n}{e_{1}, \dots , e_{r+1}} _{t} \, = \, 
t^{\nu (d)} \cdot \binom{n}{e_{1}, \dots , e_{r+1}} _{1 / t} \, = \, 
\sum _{k \in {\Bbb N}} I_{n}(d; \nu (d) - k) t^{k} \]
and hence the second part of assertion (v) readily follows. From the last remark preceding Theorem 3 it follows that
the coefficient of $t^{i}$ in $\binom{n}{e} _{t}$ is a positive integer for $0 \leq i \leq e (n - e)$. Since 
$\binom{n}{e_{1}, \dots , e_{r+1}} _{t}$ is a product of polynomials of the type $\binom{n}{e} _{t}$, we infer
that the coefficient of $t^{k}$ in it is also positive, {\it i.e.}, $I_{n}(d; k) \geq 1$, for $0 \leq k \leq \nu (d)$.
$\Box$ \\

\noindent{\bf Remarks:} 
\begin{enumerate} 
\item In view of the first remark preceding Theorem 4, the Netto-Knuth formula for $I_{n}(k)$ with $k \leq n$ (see [4])
is easily obtained from assertion (iii) of Theorem 4.
\item For $n \geq 2$, the sequence $\psi _{n} (i)$ need not be log-concave, {\it e.g.}, $\psi _{6} (5) = 1$, 
$\psi _{6} (7) = 2$ whereas $\psi _{6} (6) = 0$.
\item Define \[ \eta (d) := \frac{e_{1} (e_{1} - 1) + \cdots + e_{r+1} (e_{r+1} - 1)}{2} \]
Also, for a positive integer $k$, define
\[ [k]+ \, := \, \{ 0 \leq i \leq k \mid \psi _{n} (i) > 0 \} \;\;\;\; \mbox{and} \;\;\;\;
[k]- \, :=  \, \{i \in [k] \mid \psi _{n} (i) < 0 \}.  \]
From [1; Corollary 3.5], we get
\[\frac{1}{ e_{1}! \cdots e_{r+1}! } \binom{n - 1 + m}{n - 1} \, \leq \, D_{\varepsilon (d, n)}(m) \, \leq \,
\frac{1}{ e_{1}! \cdots e_{r+1}! } \binom{n - 1 + \eta (d) + m}{n - 1} \]
for all $m \in {\Bbb N}$. Hence $I_{n}(d; k)$ is bounded below by
\[ \frac{1}{\prod _{i=1}^{r+1} e_{i}!} \left  \{ \sum _{i \in [k]+} \psi _{n} (i) \binom{n - 1 + k - i}{n - 1} + 
\sum _{i \in [k]-} \psi _{n} (i) \binom{n - 1 + \eta(d) + k - i}{n - 1} \right \} \]
and bounded above by
\[  \frac{1}{\prod _{i=1}^{r+1} e_{i}!} \left \{ \sum _{i \in [k]-} \psi _{n} (i) \binom{n - 1 + k - i}{n - 1} + 
\sum _{i \in [k]+} \psi _{n} (i) \binom{n - 1 + \eta (d) + k - i}{n - 1} \right \}. \]
In particular, if $k \leq n$ (or by symmetry, $k \geq n (n - 1) /2$), then these bounds are explicit.
Observe that if $d: 1 < \cdots < n - 1$, or equivalently if $\eta(d) = 0$, then each of the two bounds coincides 
with $I_{n}(k)$. Suppose $n \geq 3$, $k \geq 2$ and $d$ is such that $\eta(d) \geq 1$. Then, since $\psi _{n}(1) = -1$
and   
\[ I_{n}(k) \, \leq \, \frac{n (n - 1)}{2} \, \leq \, \frac{n - 1}{k}\cdot \binom{n - 1 + k - 1}{n - 1}, \]
it is easy to deduce that the above lower bound is $\leq 0$. The above upper bound can be $ < I_{n}(k)$, {\it e.g.}, 
a MAPLE computation shows that in the case of $d_{0} = 0 < d_{1} = 1 < d_{2} = 10 = n$, we have $I_{10}(12) = 47043$ 
whereas the above upper bound evaluates to a number less than $44871$ and $I_{10}(20) = 230131$ whereas the same upper 
bound evaluates to a number less than $182032$. Under refinements of $d$, the above upper bound may increase or decrease, 
{\it e.g.}, letting $n = 5$, $k = 6$ if $d: 2$, then this upper bound is $ < 84$, if $d: 1 < 2$, then it is $104$ and 
if $d: 1 < 2 < 3$, then it is $77$. 
\item For fixed $n$, the sequence $I_{n}(k)$ is known to be log-concave (since its generating function is a product
of easily verified log-concave polynomials). In contrast, for fixed $n$ and $d$ with $\eta (d) \geq 1$, the sequence 
$I_{n}(d; k)$ need not be log-concave, {\it e.g.}, for $k = 0, 1, 2, \dots \dots$, we have 
$I_{7}(2 < 4; k):\; 1, 2, 5, 8, 13, \dots \dots$. \\
\end{enumerate}

\noindent {\bf REFERENCES:}
\begin{enumerate}
\item Agnarsson, G., {\it On the Sylvester denumerants for general restricted partitions}, Proceedings of the 
Thirty-third Southeastern International Conference on Combinatorics, Graph Theory and Computing (Boca Raton, FL, 2002). 
Congr. Numer. 154 (2002), 49-60.
\item Baldoni, V.; Berline, N.; De Loera, J. A.; Dutra, B. E.; Koppe, M.; Vergne, M., {\it 
Coefficients of Sylvester's denumerant}, Integers 15 (2015), Paper No. A11, 32 pp. 
\item Knuth, D., {\it Subspaces, subsets, and partitions}, J. Combinatorial Theory 10 (1971), 178-180.
\item Knuth, D., {\it The Art of Computer Programming}, Addison-Wesley, Reading, MA, Vol. 3.
\item Stanley, R., {\it Enumerative Combinatorics, Volume 1}, Wadsworth and Brooks/Cole, 1986, Monterey, CA.
\item Sylvester, J. and Franklin, F., {\it A constructive theory of partitions, arranged in three acts, an 
interact and an exodian}, Amer. J. Math. 5 (1882), 251-330.
\end{enumerate}

\end{document}